\documentclass[11pt]{amsart}
\usepackage{amsfonts,graphics,amscd,amssymb}

\newtheorem{proposition}{Proposition}[section]
\newtheorem{theorem}[proposition]{Theorem}
\newtheorem{definition}[proposition]{Definition}

\newtheorem{lemma}[proposition]{Lemma}

\newtheorem{remark}[proposition]{Remark}

\newcommand{\la}{\langle}
\newcommand{\ra}{\rangle}

\newcommand{\bd}{\partial}
\newcommand{\x}{\times}

\newcommand{\Vol}{\operatorname{Vol}}

\newcommand{\Supp}{\operatorname{Supp}}

\newcommand{\cA}{{\mathcal A}}

\newcommand{\cC}{{\mathcal C}}

\newcommand{\QQ}{{\mathbb Q}}
\newcommand{\RR}{{\mathbb R}}
\newcommand{\TT}{{\mathbb T}}

\newcommand{\ZZ}{{\mathbb Z}}

\begin{document}
\title[Density of ergodic solenoids]{Ergodic solenoidal homology: Density of ergodic solenoids}

\author[V. Mu\~{n}oz]{Vicente Mu\~{n}oz}
\address{Instituto de Ciencias Matem\'aticas CSIC-UAM-UC3M-UCM, Serrano 113 bis, 28006 Madrid, Spain}

\address{Facultad de
Matem\'aticas, Universidad Complutense de Madrid, Plaza de Ciencias
3, 28040 Madrid, Spain}

\email{vicente.munoz@imaff.cfmac.csic.es}

\author[R. P\'{e}rez Marco]{Ricardo P\'{e}rez Marco}
\address{CNRS, LAGA UMR 7539, Universit\'e Paris XIII,  
99 Avenue J.-B. Cl\'ement, 93430-Villetaneuse, France}

\email{ricardo@math.univ-paris13.fr}

\date{5 March, 2009}



\begin{abstract}
  A measured solenoid is a laminated space endowed with a tranversal measure invariant by holonomy.
  A measured solenoid immersed in a smooth manifold produces a closed current (known as
  generalized Ruelle-Sullivan current). Uniquely ergodic solenoids are those for which
  there is a unique (up to scalars) transversal measure. It is known that for any smooth manifold,
  any real homology class is represented by a uniquely ergodic solenoid. In this paper, we prove
  that the currents associated to uniquely ergodic solenoids are dense in the space of closed currents,
  therefore proving the abundance of such objects.
\end{abstract}

\maketitle

{\hfill {\emph{Dedicated to the 100th Anniversary of S.\ M.\ Ulam}}}

\section{Introduction} \label{sec:introduction}

In \cite{MuPM}, the authors introduced the concept of measured solenoid as a
laminated space $S$ endowed with a transversal measure $\mu$. Given a measured oriented
solenoid $S_\mu$, any smooth map $f:S\to M$ into a differentiable
manifold $M$ gives rise, by integration of forms, to a closed current
 $$
 (S_\mu,f) \in \cC_k(M)\, ,
 $$
whose associated homology class
 $$
 [S_\mu,f]\in H_k(M,\RR)\,
$$
is called the generalized Ruelle-Sullivan class.

A special role in the theory is played by the solenoids which have an ergodicity property.
A uniquely ergodic solenoid is a solenoid $S$ for which there is a unique (up to multiples)
transversal measure, and this measure has support the whole of the solenoid.
In this situation, the ``geometry'' (that is, the solenoid itself) determines the
``measure''. This has a more precise significance if we note the following: for a
uniquely ergodic immersed solenoid $f:S\to M$, take one leaf $l\subset S$, and
an exhaustion $l_R$ by compact sets such
that $f(l_R)$ can be capped off by a ``small'' cap, thus defining a homology
class $[\overline{f(l_R)}]\in H_k(M,\ZZ)$. Then the normalizations
 $$
  \frac{[\overline{f(l_R)}]}{\Vol(f(l_R))} \, \longrightarrow \, [S_\mu,f] \, .
  $$
Said otherwise, any one leaf determines the whole solenoid.

The main result of \cite{MuPM} establishes that for any manifold $M$, any homology class
$a\in H_k(M,\RR)$ can be represented by a uniquely ergodic oriented solenoid $f:S\to M$, which moreover
is inmersed (but we have to allow possible transversal self-intersections if $a\cup a\neq 0$).

In this note, we want to prove that for a given homology class $a\in H_k(M,\RR)$, the currents
associated to immersed uniquely ergodic solenoids are dense in the space of closed currents
in $\cC_k(M)$ representing $a$.

\bigskip

\noindent  \textbf{Acknowledgements.} We would like to
thank Denis Sullivan for raising the question addressed in this
paper. Partially supported through Spanish MEC grant MTM2007-63582.

\section{Definitions} \label{sec:minimal}

Let us review the main concepts introduced in \cite{MuPM}.

\begin{definition}\label{def:k-solenoid} A
$k$-solenoid, where $k\geq 0$, of class $C^{r,s}$, is a compact Hausdorff space endowed with an atlas
of flow-boxes $\cA=\{ (U_i,\varphi_i)\}$,
 $$
 \varphi_i:U_i\to D^k\x K(U_i)\, ,
 $$
where $D^k$ is the $k$-dimensional open ball, and $K(U_i)\subset \RR^l$ is the transversal
set of the flow-box. The changes of charts $\varphi_{ij}=\varphi_i\circ
\varphi_j^{-1}$ are of the form
 \begin{equation}\label{eqn:change-of-charts}
 \varphi_{ij}(x,y)=(X(x,y), Y(y))\, ,
 \end{equation}
where $X(x,y)$ is of class $C^{r,s}$ and $Y(y)$ is of class $C^s$.
\end{definition}

Let $S$ be a $k$-solenoid, and $U\cong D^k \x K(U)$ be a flow-box for $S$. The sets
$L_y= D^k\x \{y\}$ are called the (local) leaves of the flow-box. A leaf $l\subset S$ of the
solenoid is a connected $k$-dimensional manifold whose intersection with any flow-box
is a collection of local leaves. The solenoid is oriented if the leaves are oriented
(in a transversally continuous way).

A transversal for $S$ is a subset $T$ which is a finite union of transversals of flow-boxes.
Given two local transversals $T_1$ and $T_2$ and
a path contained in a leaf from a point of $T_1$ to a point of $T_2$,
there is a well-defined holonomy map $h:T_1\to T_2$. The holonomy maps form a pseudo-group.

\begin{definition} \label{def:transversal-measure}
Let $S$ be a $k$-solenoid. A transversal measure $\mu=(\mu_T)$ for
$S$ associates to any local transversal $T$ a locally finite measure
$\mu_T$ supported on $T$, which are invariant by the holonomy
pseudogroup, i.e. if $h : T_1 \to T_2$ is a holonomy map, then
$h_* \mu_{T_1}= \mu_{T_2}$.
\end{definition}

We denote by $S_\mu$ a $k$-solenoid $S$ endowed with a transversal
measure $\mu=(\mu_T)$. We refer to $S_\mu$ as a measured solenoid.
Observe that for any transversal measure $\mu=(\mu_T)$ the scalar
multiple $c\, \mu=(c \, \mu_T)$, where $c>0$, is also a transversal
measure. Notice that there is no natural scalar normalization of
transversal measures.

\begin{definition} \label{def:transversal-unique-ergodicity}
Let $S$ be a $k$-solenoid. The solenoid $S$ is
uniquely ergodic if it has a unique (up to scalars)
transversal measure $\mu$ and its support is the whole of $S$.
\end{definition}

\bigskip

Now let $M$ be a smooth manifold of dimension $n$. An immersion of a $k$-solenoid
$S$ into $M$, with $k<n$, is a smooth map $f:S\to M$ such that the differential
restricted to the tangent spaces of leaves has rank $k$ at every
point of $S$. The solenoid $f:S\to M$ is
transversally immersed if for any flow-box $U\subset S$ and chart
$V\subset M$, the map $f:U= D^k\x K(U) \to V\subset \RR^n$ is
an embedding, and the images of
the leaves intersect transversally in $M$. If moreover $f$ is injective, then
we say that the solenoid is embedded.

Note that under a transversal immersion, resp.\ an embedding,
$f:S\to M$, the images of the leaves are immersed, resp.\
injectively immersed, submanifolds.

\bigskip

Denote by
 $$
 \cC_k(M)
 $$
the space of compactly supported currents of dimension $k$
on $M$. We have the following.

\begin{definition}\label{def:Ruelle-Sullivan}
Let $S_\mu$ be an oriented measured $k$-solenoid. An immersion
$f:S\to M$
defines a generalized Ruelle-Sullivan current $(S_\mu,f)\in \cC_k(M)$ as follows.
Let $S=\bigcup_i S_i$ be a measurable partition such that each
$S_i$ is contained in a flow-box $U_i$. For $\omega\in \Omega^k(M)$, we define
 $$
 \la (S_\mu,f),\omega \ra=\sum_i \int_{K(U_i)} \left ( \int_{L_y\cap S_i}
 f^* \omega \right ) \ d\mu_{K(U_i)} (y) \, ,
 $$
where $L_y$ denotes the horizontal disk of the flow-box.
\end{definition}

In \cite{MuPM} it is proved that $(S_\mu,f)$ is a closed current. Therefore, it defines
a real homology class
 $$
  [S_\mu, f]\in H_k(M,\RR)\, .
  $$
In their original article \cite{RS}, Ruelle and Sullivan defined
this notion for the restricted class of solenoids embedded in $M$.

In \cite{MuPM}, it is proved that if $(S_\mu,f)$ is an embedded
solenoid and the transversal measure has no atoms (for instance, if $S$ has
no compact leaves), then
 $$
 [S_\mu,f] \cup [S_\mu,f]=0 \, .
 $$
So if $a\in H_k(M,\RR)$ is a homology class with $a\cup a\neq 0$, then it cannot
be represented by an embedded solenoid.

\medskip

Now introduce a Riemannian metric on $M$.
Let $(S_\mu,f)$ be a uniquely ergodic immersed oriented solenoid. The leaf-wise
volume of $S$ (induced by the metric on $M$) together with the transversal measure,
give a finite measure supported on $S$. We normalize it to have total mass $1$.
Note that this produces a unique transversal measure $\mu$.

Let $l\subset S$
be a leaf of $S$. Suppose that there is an exhaustion $l_R$ of $l$ such that
$f(l_R)$ has a small cap, that is, there is an oriented submanifold $C_R$ of dimension $k$
with boundary $\bd C_R=-\bd f(l_R)$, and satisfying
 \begin{equation} \label{eqn:llld}
  \frac{\Vol (C_R)}{\Vol f(l_R)} \, \longrightarrow \, 0\, .
 \end{equation}
Then consider the integration current  $(l_R,f)\in \cC_k(M)$, defined by
$\la (l_R,f),\omega\ra= \int_{l_R} f^*\omega$. Then we have the following:

\begin{proposition} \label{prop:limit}
 In the situation above, $(l_R,f)/{\Vol f(l_R)} \to (S_\mu,f)$.
\end{proposition}

\begin{proof}
  First note that the currents $(l_R,f)/{\Vol f(l_R)}$ are bounded. Therefore
  they accumulate. To see that the sequence converges, it is enough to check
  that there is only one accumulation point. So assume that $(l_R,f)/{\Vol f(l_R)}$
  converges, say to some current $T$.

  First note that the current $T$ is closed:
   \begin{eqnarray*}
   \la \bd T,\omega \ra &=& \lim_{R\to \infty} \frac1{\Vol f(l_R)}
   \la \bd (l_R,f),\omega \ra = \lim_{R\to \infty} \frac1{\Vol f(l_R)}
   \la \bd C_R , \omega \ra \\ &=& \lim_{R\to \infty} \frac1{\Vol f(l_R)} \int_{C_R} d\omega =0\, ,
   \end{eqnarray*}
 by (\ref{eqn:llld}).

   Clearly, $T$ is supported on $S$. By unique ergodicity, $T$ should be a multiple
   of $(S_\mu,f)$. Let us see this: $T$ is a current defined by a transversal measure
   if and only if it is a daval measure as defined in \cite{MuPM} (that is, locally it
   is the product of the Riemannian volume along leaves with a transversal measure).
   This is equivalent to $T$ being invariant by the group $G_S^0$ consisting of diffeomorphisms
   of the solenoid isotopic to the identity. So we have to check that for a leaf-wise
   vector field $X$ defined on the solenoid, it is $L_XT=0$. But for any $(k+1)$-form $\beta$ on $S$,
   $i_X\beta|_S=0$, so $\int_{l_R} i_X\beta=0$. This means that
   $\la i_X (l_R,f), \beta\ra =\la (l_R,f), i_X\beta\ra =0$. Therefore $i_XT=0$, and hence
   $$
   L_XT=di_XT+i_XdT=0\, ,
   $$
   as claimed.

  The argument above yields that $T=\lambda (S_\mu,f)$, for some $\lambda\in \RR$.
  Clearly $\lambda\geq 0$, since both currents define the same orientation at a fixed
  point of the leaf $l$.
  Consider the norm dual to the $C^0$-norm on $\cC_k(M)$. It is easy to see
  that $||(l_R,f)/{\Vol f(l_R)}||=1$. Therefore $||T||=1$. But our normalization
  implies that $||(S_\mu,f)||=1$. So $T=(S_\mu,f)$.
\end{proof}

\section{Realization theorem}\label{sec:k-solenoids}

Let $M$ be a smooth compact oriented Riemannian $C^\infty$ manifold
and let $a\in H_k(M,\RR)$ be a non-zero real $k$-homology class. The main
result of \cite{MuPM} is the construction of a uniquely ergodic oriented transversally
immersed $k$-solenoid representing a positive multiple of $a$.

The construction is as follows.
Take a collection
$C_1,\ldots, C_{b_k}\in H_k(M,{\ZZ})$ which are a
basis of $H_k(M,{\QQ})$ and such that $C_i$ are represented by immersed
submanifolds $S_i\subset M$ with trivial normal bundle and such that
all
intersections are transversal.
After switching the orientations of $C_i$ if necessary,
reordering the cycles and multiplying $a$ by a suitable positive
real number, we may suppose that
 $$
 a=\lambda_1 C_1 +\ldots +\lambda_r C_r,
 $$
for some $r\geq 1$, $\lambda_i>0$, $1\leq i\leq r$, and $\sum
\lambda_i=1$.

Let
$h:\TT\to \TT$ be a diffeomorphism of the circle which is a Denjoy
counter-example with an irrational rotation number and of class
$C^{2-\epsilon}$, for some $\epsilon>0$. Hence $h$ is uniquely
ergodic. Let $\mu_K$ be the unique invariant probability measure,
which has support on a Cantor set $K\subset \TT$. Partition the
Cantor set $K$ into $r$ disjoint compact subsets $K_1,\ldots, K_r$
in cyclic order, each of which with $\mu_K(K_i)=\lambda_i$.

Now construct the immersed solenoid $f:S\to M$ as follows. The central
part (or core) of $S$ lies inside a ball in $M$, and it is of the form
$C= S^{k-1}\x [-1,1]\x K$ embedded as
 $$
 (x,t,y) \mapsto (x,t, h_t(y))\, ,
 $$
where $h_t$ is an isotopy of $\TT$ from the identity to $h$.

\begin{figure}[h]\label{figure4}
\centering
\resizebox{10cm}{!}{\includegraphics{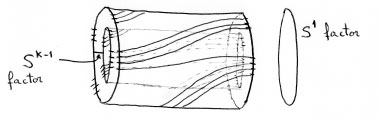}}    
\caption{The manifold $X$.}%
\end{figure}

The rest of the solenoid consists of $(S_i-(D_i^+ \cup D_i^-))\times K_i$
for $i=1,\ldots, r$, where $S_i-(D_i^+ \cup D_i^-)$ is $S_i$ with two
balls removed, and we are taking parallel copies of $S_i$ as leaves.
We glue the boundaries $\bd D_i^\pm \times K_i$ to the boundaries
$S^{k-1}\times \{ \pm 1\} \x K_i\subset \bd C$.

This gives an oriented $k$-solenoid of class $C^{\infty,2-\epsilon}$. That is, the
changes of charts \eqref{eqn:change-of-charts} are $C^{\infty}$ in $x$ and $C^{2-\epsilon}$ in $y$.

We have a global transversal $T=\{p\}\x K \subset S^{k-1}\x K\subset S$.
Identifying $T\cong K$, the holonomy pseudo-group is generated by
$h:K\to K$. Hence $S$ is uniquely ergodic. Let $\mu$ denote the
tranversal measure corresponding to $\mu_K$.

\begin{figure}[h]\label{figure4}
\centering
\resizebox{10cm}{!}{\includegraphics{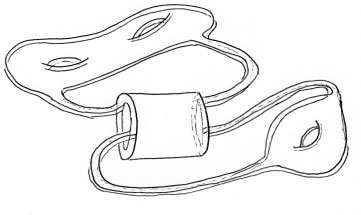}}    
\caption{The solenoid $S$.}
\end{figure}

\medskip

It remains to prove that $[S_\mu,f]=a$ (up to an scalar multiple). This is done
by coupling against suitable test forms, as done in \cite{MuPM}.

\section{Density of uniquely ergodic solenoids} \label{sec:density}

Now we move on to the main result of this paper. We consider the space of
currents $\cC_k(M)$ with the weak topology (that is, the topology as dual
space to $\Omega^k(M)$ with the Fr\'echet topology given by the $C^\infty$-convergence
of forms).

\begin{theorem} \label{thm:main}
 Let $T\in \cC_k(M)$ be a (compactly supported)
 closed current representing a homology class $a\in H_k(M,\RR)$.
 Then there is an oriented immersed
 uniquely ergodic solenoid $f:S\to M$ such that the current $(S_\mu,f)$ is
 as close to $T$ as we want, and $[S_\mu,f]=a$.
\end{theorem}

\begin{proof}
We start by considering the oriented, transversally immersed,
measured solenoid $f_a:S_a\to M$, constructed in Section \ref{sec:k-solenoids}.
It satisfies that $[S_a,f_a]=a$. Then
 $$
 T-(S_a,f_a) \in \cC_k(M)
 $$
is an exact current. Therefore there is some $T'\in \cC_{k+1}(M)$ with
 $$
 T-(S_a,f_a)= \bd T'\, .
 $$
By a regularization
procedure (or using the density of the forms in the space of currents),
we can find a smooth $(n-k-1)$-form $\beta$ (thought of as a current, by duality)
such that $T'-\beta$ is very small.
So
  $$
  T \approx (S_a,f_a) + d\beta\, .
  $$

Now take an open cover $\{U_i\}$ of $M$ by coordinate open sets, and a refinement
$\{V_i\}$ such that $\overline{V}_i\subset U_i$. Consider a partition of
unity $\rho_i$ subordinated to $V_i$, and let $\beta_i=\rho_i\beta$.
Writing in coordinates $(x_1,\ldots,x_n)$ for $U_i$, we have
  $$
  \beta_i= \sum_{i_1<\ldots< i_{n-k-1}} h_{i_1\ldots i_{n-k-1}} \,
  dx_{i_1}\wedge \ldots \wedge dx_{i_{n-k-1}} \, .
  $$
Consider now a bump function $\tilde\rho_i$,
which is one on $V_i$ and zero off $U_i$.
Then
  $$
  \beta_i= \sum_{i_1<\ldots< i_{n-k-1}} \tilde\rho_i\, h_{i_1\ldots i_{n-k-1}}
  \,  d(\tilde\rho_i\, x_{i_1})\wedge \ldots
  \wedge d(\tilde\rho_i\, x_{i_{n-k-1}})\, ,
  $$
and
  $$
  d\beta_i= \sum_{i_1<\ldots< i_{n-k-1}} d(\tilde\rho_i\, h_{i_1\ldots i_{n-k-1}})
  \wedge  d(\tilde\rho_i\, x_{i_1})\wedge \ldots
  \wedge d(\tilde\rho_i\, x_{i_{n-k-1}})\, .
  $$

As $d\beta=\sum d\beta_i$, we have that $d\beta$ can be written as a sum
of terms of the form
 \begin{equation} \label{eqn:aqa}
 df_1^i\wedge \ldots \wedge df_{n-k}^i\, ,
 \end{equation}
where $f_1^i,\ldots, f_{n-k}^i$ are (compactly supported)
globally defined functions on $M$.

If $f_1,\ldots, f_{n-k}$
are compactly supported functions on $M$, write
 \begin{equation}\label{eqn:F}
 F=(f_1,\ldots, f_{n-k}):M\to \RR^{n-k} \, .
 \end{equation}
 Then the
 form $\alpha=df_1\wedge \ldots \wedge df_{n-k}= F^* \omega$, where
 $\omega=dx_1\wedge \ldots \wedge dx_{n-k}$ is the volume form on $\RR^{n-k}$.

\begin{lemma}\label{lem:alpha}
 If $\alpha=df_1\wedge \ldots \wedge df_{n-k}$, with $F$ as in \eqref{eqn:F}, then
 there is an embedded solenoid with trivial holonomy and transversal
 Cantor sets whose associated current is as close to $\alpha$ as we want.
\end{lemma}

\begin{proof}
  Consider the set $C=\{x\in M\, ; \, \alpha(x)=0\}$. Then $F(C)\subset \RR^{n-k}$ is the set
  of critical values of $F$, which has zero measure by Sard's theorem. Denote by $K=F(M)\subset
  \RR^{n-k}$, and note that it is compact since $F$ is compactly supported. Now let $U$ be
  a small open neighbourhood of $F(C)$. Let us see that
   \begin{equation}\label{eqn:a}
   F^*\omega - F^*(\omega|_{K-U})= F^*(\omega|_U)
     \end{equation}
  is small (as a current in $M$). Fix $\epsilon>0$, and consider $C_\epsilon=
  \{x\in M \,  ; \, |\alpha(x)|<\epsilon\}$.
For $\beta \in \Omega^k(M)$, we have
   \begin{equation}\label{eqn:1}
   \left|\int_{F^{-1}(U) \cap C_\epsilon} F^*\omega \wedge \beta \right| \leq
   \epsilon\, \Vol(\Supp(F)) |\beta|\, ,
   \end{equation}
  so it is small. In the complement $M-C_\epsilon$, the leaves
  $F^{-1}(x)$ have uniformly bounded volume. This follows from the fact
  that the $dF$ has norm $\geq \epsilon$ in the normal directions to the leaves,
  and the total volume of $M-C_\epsilon$ is bounded (recall that $M$ may be
  non-compact, but as $F$ is compactly supported, $\overline{M-C}$ is compact,
  and hence $M-C_\epsilon$ is also compact). Therefore
   \begin{equation}\label{eqn:2}
   \left|\int_{F^{-1}(U) \cap (M-C_\epsilon)} F^*\omega \wedge \beta \right| =
   \left| \int_{U} \left(\int_{F^{-1}(x)} \beta \right)
   dx_1\wedge \ldots dx_{n-k} \right| \leq C \Vol(U) \, ,
     \end{equation}
  which is small by taking $U\supset C$ very small. Adding up \eqref{eqn:1} and
  \eqref{eqn:2}, we get that \eqref{eqn:a} is small.

The current $F^*(\omega|_{K-U})$ is actually a solenoid. The leaves
$F^{-1}(x)$ are compact submanifolds, and the transversal
measure is the pull-back of $\omega$ to the transversals. Note
that $K-U$ is a global transversal, and
that the holonomy is trivial.
Now consider a measure $\mu$, supported on a Cantor set which approximates $\omega|_{K-U}$.
Then
 $$
 \la F^*\mu,\beta \ra =\int_{K-U} \left( \int_{F^{-1}(y)} \beta \right) d\mu(y)
 = \la \mu, F_*\beta\ra\, .
$$
Since the norm of $dF$ is bigger than a fixed constant,
the norm of $F_*\beta$ is bounded by a constant
times the norm of $\beta$. Therefore
 $$
 |\la F^*\mu,\beta \ra - \la F^*\omega,\beta\ra |
 = |\la \mu-\omega, F_*\beta \ra |\leq \varepsilon |F_*\beta|
\leq C\varepsilon |\beta|\, ,
 $$
for some small $\varepsilon>0$. This completes the proof of the lemma.
\end{proof}

Using this lemma, we get a collection of solenoids, with transversals being
Cantor sets, and whose union defines the sought current. All but one have
trivial holonomy, and the remaining one is uniquely ergodic. Our next task
is to do a surgery to get a connected (and uniquely ergodic) solenoid.
First note that we may decompose a solenoid with trivial holonomy in
smaller chunks so that its total (transversal) mass is smaller than that
of the solenoid $S_a$.

\begin{proposition} \label{prop:surgery}
 Let $f_i: S_{i,\mu_i} \to M$, $i=1,2$, be two 
 immersed oriented
 measured solenoids of the same dimension $k$, and with Cantor transversal
 structure. Suppose that there are two
disjoint discs $D^{n-k}_i\subset M$ such that $T_i= D^{n-k}_i\cap S_i$
is a global transversal for $S_i$. Let $h_i$
be the holonomy of $S_i$.
Suppose also that there is a diffeomorphism $\varphi: T_1\to T_2$
which preserves the measures. Then there is an 
immersed oriented
measured solenoid $f:S_\mu \to M$ with
global transversal $T=T_1$, and holonomy $h=\varphi^{-1}\circ h_2
\circ \varphi \circ h_1$, such that
the associated currents satisfy
 $$
 (S_\mu,f) \approx (S_{1,\mu_1},f_1)+ (S_{2,\mu_2},f_2) \, .
 $$
\end{proposition}

\begin{proof}
  Take flow boxes $D^k\times T_i$ for $S_i$, with $\{0\}\x T_i$
  corresponding to $T_i$, remove the interiors
  $D^k_{1/2} \times T_i$, and glue the boundaries, to get the
  required (abstract) solenoid
  $$
  S= (S_1 -D^k_{1/2} \times T_1) \cup_{\bd D^k_{1/2} \times T_1\cong_\varphi \,
  \bd D^k_{1/2} \times T_2} (S_2- D^k_{1/2} \times T_2) \, .
  $$
  The transversal measure is induced by $\mu_1$ (or $\mu_2$).

It remains to define the immersion $f:S\to M$. On $S_i-D^k_1 \times T_i$, it is defined
to be equal to $f_i$. Now fix $y_0\in T_1$ and $y_0'=\varphi(y_0)\in T_2$, and consider
a path $\gamma:[0,1]\to M$ from $y_0$ to $y_0'$ transversal to the leaves through the
end-points. 
Fatten up $\gamma$ to a map $\hat\gamma:
D^k \x [0,1] \to M$, which matches the leaves at the end-points. We can take a small
open set $U\subset D_1^{n-k}$ and the corresponding $\varphi(U)\subset D_2^{n-k}$, so that we can
extend $\hat\gamma$ to a map
 $$
 \tilde\gamma: D^k\x [0,1] \x (U\cap T_1) \to M\, .
 $$
Repeating this for a finite cover of $T_1$, we get an open subset $V\supset T_1$, and a map
 $$
 \tilde\gamma: D^k \x [0,1] \x T_1 \to M\, .
 $$

Now define $f:S\to M$ as follows. Firstly, send $(D_1^k-D_{3/4}^k)\x T_i$ to
$f_i((D^k-D_\epsilon^k)\x T_i)$, and secondly, define $f$ on
$((D_{3/4}^K - D_{1/2}^k) \x T_1) \cup_\varphi ((D_{3/4}^k- D_{1/2}^k)
\x T_2)$ as $\tilde\gamma|_{S^k_\epsilon \x [0,1]\x T_1}$.
Finally, one has to smooth out the corners, but this is achieved with a very small perturbation.

To end the proof of the proposition, we need to consider the difference current
 $$
 T''= (S_\mu,f)-(S_{1,\mu_1},f_1)-(S_{2,\mu_2},f_2)\, .
 $$
The leaves $L_y''$ of this solenoidal current are diffeomorphic to the
middle portion, plus two small caps
 $$
 L_y''= \tilde\gamma (S_\epsilon^k \x [0,1]\x \{y\}) \cup f_1(D_\epsilon^k\x\{y\}) \cup
   f_2(D_\epsilon^k\x\{\varphi(y)\}), \qquad y\in T_1\, .
 $$
So
 $$
 \la T'', \beta\ra=\int_{T_1} \left(\int_{L_y''}\beta\right) d\mu(y)\, .
 $$
The volume of $L_y''$ is bounded by a constant times
$\epsilon$, so it is very small for $\epsilon$ small, as required.
\end{proof}

We end up the proof of Theorem \ref{thm:main}
as follows. We have the solenoid $(S_a,f_a)$ and a
collection of solenoids $(S_i,f_i)$ approaching \eqref{eqn:aqa}. Decompose $S_i$ into
smaller solenoids so that we can assume that the total transversal measure of $S_i$ does
not exceed that of $S_a$. We can use Proposition \ref{prop:surgery} to glue each $S_i$ (once at
a time) to $S_a$ (note that Proposition \ref{prop:surgery}
works with the transversal of $S_i$ and a sub-transversal
of $S_a$). Note here that we have to construct a diffeomorphism between the transversal Cantor
sets. For this it may be necessary to arrange the transversal Cantor set for $S_i$ in such
a way that it is diffeomorphic to the transversal Cantor set of $S_a$, but this is not
problematic, since we were only requiring that the transversal measure of $S_i$
approximates the smooth transversal measure of the foliation that $S_i$ is approximating.

The holonomy of the resulting solenoid is that of $S_a$, hence uniquely ergodic.
This completes the proof of Theorem \ref{thm:main}.
\end{proof}

\begin{remark}
 An immersed solenoid $f:S\to M$ is said to have a trapping region  if there is a
 ball $B$ in $M$ such that all holonomy phenomenon of $S$ occur inside $B$,
 and the holonomy is generated by a single map (see \cite{MuPM} for a precise definition).
 All the solenoids constructed in Sections \ref{sec:k-solenoids} and \ref{sec:density}
 have a trapping region. This is relevant since the more restricted is the class of
 solenoids that we use, the better for possible future applications.
\end{remark}

\end{document}